\documentclass[sn-mathphys,Numbered, 12pt]{sn-jnl}

\usepackage{graphicx}%
\usepackage{multirow}%
\usepackage{amsmath,amssymb,amsfonts}%
\usepackage{amsthm}%
\usepackage{mathrsfs}%
\usepackage[title]{appendix}%
\usepackage{xcolor}%
\usepackage{textcomp}%
\usepackage{manyfoot}%
\usepackage{booktabs}%
\usepackage{algorithm}%
\usepackage{algorithmicx}%
\usepackage{algpseudocode}%
\usepackage{listings}%
\usepackage{multirow}
\usepackage{graphicx,epstopdf}
\usepackage{graphicx}
\usepackage{amsfonts}
\usepackage{amsmath}
\usepackage{mathrsfs}
\usepackage{amssymb}
\usepackage[T1]{fontenc}
\usepackage{mathrsfs}
\usepackage{dsfont}
\usepackage{xcolor}
\usepackage{float}
\usepackage{caption}
\usepackage{geometry}
\newtheorem{lemma}{Lemma}
\usepackage{graphicx}
\usepackage{mathptmx}      
\usepackage{array}

\theoremstyle{thmstyleone}%
\newtheorem{theorem}{Theorem}
\newtheorem{proposition}[theorem]{Proposition}%

\theoremstyle{thmstyletwo}%
\newtheorem{remark}{Remark}%

\theoremstyle{thmstylethree}%

\begin{document}

\title[]{On kernel mode estimation under RLT and WOD model}


\author*[1]{\fnm{Mohamed Kaber} \sur{El Alem}}\email{med.kaber.nema@gmail.com}

\author[1]{\fnm{Zohra } \sur{Guessoum}}\email{zguessoum@usthb.dz}
\equalcont{These authors contributed equally to this work.}

\author[1]{\fnm{Abdelkader} \sur{ Tatachak}}\email{atatachak@usthb.dz}
\equalcont{These authors contributed equally to this work.}

\affil[1]{\orgdiv{Laboratory MSTD, Department of Probability and Statistics}, \orgname{University of Sciences and Technology Houari  Boumediene (USTHB)}, \orgaddress{\street{BP 32}, \city{El-Alia}, \postcode{16111}, \state{Algiers}, \country{Algeria}}}




\abstract{Let $(X_N)_{N\geq 1}$ denote a sequence of real random variables and let $\vartheta$ be the mode of the random variable of interest $X$.
In this paper, we study the kernel mode estimator (say) $\vartheta_n$ when the data are widely orthant dependent (WOD) and subject to Random Left Truncation (RLT) mechanism. We establish the uniform consistency rate of the density estimator (say) $f_n$ of the underlying density $f$ as well as the almost sure convergence rate of $\vartheta_n$. The performance of the estimators are illustrated via some simulation studies and applied on a real dataset of car brake pads.  
}
\maketitle
\section{Introduction and motivation}
In statistical inference, there are some situations where we are unaware of the distribution governing the population, particularly when it does not align with any parametrized family of laws. This situation often better reflects the complexities of reality. Hence, in such situations, it becomes challenging to estimate one or more parameters defining the underlying law. Instead, it is imperative to estimate the target law through its density (say) $f$, employing methods such as non-parametric estimation.
Parzen (\cite{P}, 1962), Rosenblatt (\cite{R}, 1956), and Silverman (\cite{Silv}, 1986) presented early results on kernel density estimation, and subsequent research has further explored this area.\\
Density function offers several advantages; notably, estimating $f$ provides a meaningful approach to estimating various characteristics, including mean, variance, moments, quartiles, etc. This facilitates to visualize the underlying distribution function, enabling the identification of areas with high or low probabilities. In recent years, there has been considerable interest in density estimation. One of the most widely employed techniques is the kernel method, pioneered by Rosenblatt (\cite{R}, 1956) and Parzen (\cite{P}, 1962). Their groundbreaking work introduced a class of estimates entirely determined by a kernel function $K$ and a smoothing parameter $h_n$.\\
The mode is one of the measures of central tendency used to identify the most frequently occurring values. For a probability density function $f$, it is the value at which $f$ attains a maximum and to express its importance as a robust parameter, Bickel (\cite{Bickel}, 2002) concluded that, while the median is resistant to outliers, the mode is immune to them; also, it is a safer measure of location when the data may suffer from the latter.
The problem of mode estimation may be considered as a direct consequence of density estimation and has been extensively addressed in the literature. The most popular mode estimator is the well-known  one proposed and studied by Parzen (\cite{P}, 1962).\\
In many situations it may not be possible to observe the data completely and we do not have sufficient information about the individuals before the time of data recruitment. Among the different forms in which incomplete data appear, censoring and truncation are two common forms that are practically involved in survival analysis and reliability theory. In particular, we will focus on the case of  RLT model. This type of data appears in medical studies, mainly in the analysis of the life span of patients with a particular disease, it also occurs in industrial and insurance studies.  Woodroofe (\cite{W}, 1985) reviewed examples from astronomy and economics where such data can occur.\\
The problem of estimating the unconditional/conditional mode of a probability density has been addressed in statistics literature, and number of recent papers dealt with this topic. To quote a few of them. A kernel estimation procedure is used in most of these works. We can refer to Woodroofe (\cite{W}, 1985) and Stute (\cite{S}, 1993) where the distribution of left-truncated data was estimated and the asymptotic properties of the estimator were derived.\\
Recall that under RLT model in both iid, $\alpha$-mixing and associated hypotheses, Ould Sa\"{\i}d and Tatachak (\cite{OT1,OT}, 2009) and Guessoum and Tatachak (\cite{GT}, 2020) established strong consistency rates for kernel mode estimators, while Ferrani et al. (\cite{Yacine}, 2016) studied the strong uniform convergence of the kernel density and mode estimate for associated and censored data. The asymptotic normality of the kernel mode estimator under RLT and strong mixing condition was studied by Benrabah et al. (\cite{Benrabah}, 2015).\\  
An essential inquiry revolves around whether the consistency property of the proposed kernel mode estimator is preserved when dealing with truncated WOD data. Addressing this question constitutes the primary objective of our research.\\
On the other hand, it is known that in statistical applications the independence assumption is not always reasonable. This is why various dependent structures have been introduced in the last decades such as negatively associated (NA), negatively superadditive dependent (NSD), negatively orthant dependent (NOD), extended negatively dependent (END). 
Dependence relations between random variables are one of the most studied topics in statistics, such as strong mixing, association and WOD conditions. 
one of the new dependent structures that has attracted the interest of statisticians has been named WOD structure of random variables, which contains most negative dependent random variables, some positive dependent random variables and other random variables. It is also useful for the search of ruin model in the field of probability. In the literature, it has been pointed out that NA implies NSD, NSD implies NOD, NOD implies END, and END implies WOD and that the reverse is generally not true. For more details, we refer readers to Joag-Dev and Proschan (\cite{JP}, 1983), Hu (\cite{Hu}, 2000), Lehmann (\cite{Lehman}, 1966) and Liu (\cite{Liu2009}, 2009). This new dependency structure was introduced by Wang et al. (\cite{Wang2013}, 2013). \\																				
For a sequence of random variables $\{X_N, N\geq 1\}$, if there exists a finite real sequence $\{g_U(N), N\geq 1\}$ satisfying for $N\geq 1$ and $\forall x_i \in (-\infty, +\infty)$, $1\leq i \leq N$,
$$
\mathbb{P}(X_1>x_1,X_2>x_2,\cdots,X_N>x_N)\leqslant g_U(N)\prod_{i=1}^N \mathbb{P}(X_i>x_i)
$$ 
and if there also exists a sequence $\{g_L(N), N\geq 1\}$ satisfying for $N\geq 1$ and $\forall x_i \in (-\infty, +\infty)$, $1\leq i \leq N$,
$$
\mathbb{P}(X_1\leqslant x_1,X_2\leqslant x_2,\cdots,X_N\leqslant x_N)\leqslant g_L(N)\prod_{i=1}^N \mathbb{P}(X_i\leqslant x_i)
$$ 
Then the random sequence $\{X_N, N\geq 1\}$ is called Widely Orthant Dependent (WOD) with the dominating coefficients $g(N)=\max \{g_U(N), g_L(N)\}$, where $g_U(N) \geq 1$ and $g_L(N) \geq 1$.\\
\noindent
We can also refer to some work on non-parametric estimation of the density function based on WOD samples. For example, Shi and Wu (\cite{Shi}, 2014) studied the strong consistency of the kernel density estimator  for identically distributed WOD samples. Li et al. (\cite{Li}, 2015) studied the strong pointwise consistency of a type of recursive kernel estimator for WOD samples. Recently, Wang et al. (\cite{WANG}, 2022) established the convergence rate of the kernel density estimator for widely orthant dependent random variables. To our knowledge, no results exist on the nonparametric density estimator for incomplete and WOD data, except the results by  Wu et al. (\cite{CWOD}, 2024) stated in a censoring and WOD context. This work aims to extend previous results to truncated and WOD data.\\			
The goal of this study is to investigate the asymptotic behaviors of the kernel estimator of the density function $f$, particularly in scenarios where the data are subject to random left truncation and exhibit a dependence structure known as WOD. As an application, we present the strong uniform rate of the simple estimate of the mode. \\
It is noteworthy that the random left truncation (RLT) mechanism preserves the WOD property. If the original sequence of interest $\{(X_i,Y_i);\, i = 1,\cdots, N\}$ is WOD, then the observed sequence $\{(X_i,Y_i);\, i = 1,\cdots, n\}$ (where $n \leqslant N$) is also WOD. In other words, any subset of WOD random variables remains WOD, and this property follows from the definition of WOD random variables.\\
The organization of this paper is as follows. In section 2, the necessary notations are introduced and some preliminaries are listed. In section 3, the main asymptotic results are presented. In section 4, we perform a simulation study. In section 5, the main results of section 3 are proved.
\section{Preliminaries and notations}
Let $\{X_i :1\leq i\leq N \}$ be a sequence of real survival times in a life table defined on a common probability space $(\Omega,\mathcal{A},\mathbb{P})$.These random variables are not assumed to be mutually independent; instead, they have a continuous but unknown common marginal distribution function (df) $F$ and marginal density $f$. Let $\{Y_i :1\leq i\leq N\}$ be a sequence of truncating random variables with a common continuous and unknown marginal df $G$. In addition, the $Y_i's, 1\leq i\leq N$ are assumed to be independent of the $X_i's, 1\leq i\leq N$. In the RLT  model, the pair $(X, Y)$ is observed if $X\geq Y$, otherwise we have no information about them. Thus, among the $N$ random variables, we can only observe those $n=\sum_{i=1}^{N}\mathds{1}_{X_i\geq Y_i}$ pairs $(X_k,Y_k); 1\leq k\leq n$.  Without confusion, we will denote by $\{(X_i, Y_i ): 1\leq i\leq n \}$ the observed  pairs.
 The size of the actually observed sample, $n$, is a random variable. Define $\alpha := P(X \geq Y )$, it is clear that if $\alpha = 0$, no data can be observed and so throughout this paper we assume that $\alpha > 0$.\\
 In the rest of this paper, our results will not be stated with respect to the probability measure $\mathbb{P}$ (related to sample $N$) but with respect to the probability measure $\mathbf{P}$ (related to sample $n$). Similarly, $\mathbb{E}$ and $\mathbf{E}$ denote the expectation operators related to $\mathbb{P}$ and $\mathbf{P}$, respectively. In the framework of the left truncation model, the joint conditional distribution of an observation $(X, Y)$, becomes
\begin{eqnarray}
     \mathbf{H}(x,y) &=& \mathbf{P}(X\leqslant x,Y\leqslant y ) \notag \\
              &=& \mathbb{P}(X\leqslant x,Y\leqslant y | X\geqslant Y ) = \frac{1}{\alpha}\int_{-\infty}^{x}G(y\wedge z)dF(z), \notag 
    \end{eqnarray}
where \; $y\wedge z = \min(y,z)$. The marginal conditional distributions are defined by
 $$ \mathbf{F}(x):=\frac{1}{\alpha}\int_{-\infty}^{x}G(z)dF(z) \quad \text{and} \quad 
  \mathbf{G}(y):=\frac{1}{\alpha}\int_{-\infty}^{+\infty}G(y\wedge z)dF(z) $$
So the marginal conditional probability density function of $X$ is 
\begin{equation}\label{equat1}
d\mathbf{F}(x)=\frac{1}{\alpha}G(x)f(x).
\end{equation}  
As it is discussed before, we are interested in estimating $f(\cdot)$, so from (\ref{equat1}) we have
\begin{equation}\label{equat2}
f(x)=\frac{\alpha}{G(x)}d\mathbf{F}(x),
\end{equation}
and we use the following kernel density estimator for $f(x)$, deduced from (\ref{equat2})
\begin{equation}
f_n(x)=\frac{\alpha}{nh_n}\sum_{i=1}^{n}K\left(\frac{x-X_i}{h_n} \right)\frac{1}{G(X_i)},
\end{equation}
where $\{h_n\}_{n\geq 1}$ is a bandwidth sequence, such that $h_n\rightarrow 0$ as $n\rightarrow \infty$, and $K(\cdot)$ is some kernel function.\\
When $G$ is known, $f_n(x)$ can be used to estimate the common density of the interest variables. However, in most practical cases $G$ is unknown and can be remplaced by the Lynden-Bell estimator $G_n(\cdot)$.\\
Let $ a_w=\inf\{u: W(u)>0\} \quad \text{and} \quad  b_w=\sup\{u: W(u)<1\}$ be  the lower and upper bounds of the $W$ distribution function support. As in Woodroofe (\cite{W}, 1985), $F$ and $G$ can be estimated completely only if
$$ a_G\leqslant a_F,\qquad b_G\leqslant b_F \quad \text{and} \quad \int_{a_F}^{\infty}\frac{dF}{G}<\infty $$
Let $C(\cdot)$ be the function  defined by
\begin{equation}\label{eq3}
C(x):= \mathbf{P}(Y \leqslant x \leqslant X )=\mathbf{G}(x)-\mathbf{F}(x)=\frac{1}{\alpha}G(x)[1-F(x)].
\end{equation}
The functions $\mathbf{F}$, $\mathbf{G}$ and $C$ can be estimated empirically by 
$$ \mathbf{F_{n}}(x)=\frac{1}{n}\sum_{i=1}^{n}\mathds{1}_{\{X_i \leqslant x \}}, \quad  \quad \mathbf{G_{n}}(y)=\frac{1}{n}\sum_{i=1}^{n}\mathds{1}_{\{Y_i \leqslant y\}} \qquad \text{and} \qquad
 C_n(x)=\frac{1}{n}\sum_{i=1}^{n}\mathds{1}_{\{Y_i \leqslant x \leqslant X_i\}} $$ respectively, where $\mathds{1}_A$ designates the indicator function of the  set $A$. \\
Lynden-Bell (\cite{L}, 1971) constructed a nonparametric estimators of F and G given by 
\begin{equation}\label{eq4}
F_n(x)=1-\prod_{i:X_i\leqslant x}\left[\frac{nC_n(X_i)-1}{nC_n(X_i)}\right] \quad \text{and} \quad G_n(y)=\prod_{i:Y_i>y}\left[\frac{nC_n(Y_i)-1}{nC_n(Y_i)}\right].
\end{equation}
According to (\ref{eq3}) and replacing F and G by their respective non-parametric maximum likelihood estimator, we can consider the estimator of $\alpha$, namely
\begin{equation}\label{eq5}
\alpha_n(x) = \frac{G_n(x)[1-F_n(x)]}{C_n(x)}\mathds{1}_{\{C_n(x)\neq 0\}}=:\alpha_n.
\end{equation}
He and Yang (\cite{HY}, 1998), proved that $\alpha$ does not depend on $x$ and they have shown that it is strongly consistent for $\alpha$.\\
According to (\ref{eq3}) and (\ref{eq5}), we are now in a place to present a more applicable estimator of $f$, noted $\hat{f}_n$ and  defined by
\begin{equation}\label{e1}
\hat{f}_n(x)=\frac{\alpha_n}{nh_n}\sum_{i=1}^{n}K\left(\frac{x-X_i}{h_n} \right)\frac{1}{G_n(X_i)}\mathds{1}_{\{G_n(X_i)\neq 0\}}.
\end{equation}
Assume now that $f$ is unimodal and denote by $\vartheta$ its mode, which is defined by the following equation
$$
\vartheta = \arg  \max_{x\in \mathbb{R}} f(x).
$$
The kernel estimator of $\vartheta$ is defined as the random variable $\vartheta_n$ that maximizes the kernel estimator $\hat{f}_n(x)$ of $f(x)$, i.e.
\begin{equation}\label{e2}
\vartheta_n = \arg  \max_{x\in \mathbb{R}} \hat{f}_n(x).
\end{equation}
Recall that asymptotic results for (\ref{e1}) and (\ref{e2}), in both independent and identically distributed (iid) and strong mixing condition cases have been stated in (Ould~Sa\"{i}d and Tatachak  (\cite{OT1}, \cite{OT}, 2009), Benrabah et~al.  (\cite{Benrabah}, 2015)) under random left truncation. In association condition case Guessoum and Tatachak (\cite{GT}, 2020) establish the strong uniform consistency with a rate of a kernel function estimator, in (\ref{e1}), when the variable of interest is subject to random left truncation. The results obtained in this paper extend these authors' results to a more general dependency structure known as WOD.
\section{Theoretical results}
In this section we will present our main results and provide some necessary assumptions, which will be used to establish these latter. Before stating our results, we need a few preliminary elements. In the following, all limit relations are expressed in $n \rightarrow \infty$. For two positive functions
$u(n)$ and $v(n)$, we note  $u(n) = O(v(n))$ if $\limsup\, u(n)/v(n) < \infty$. Furthermore, we consider a compact $D :=[a, b]$ such that $a_G \leqslant a_F < a < b < b_F$. Due to these restrictions and without loss of generality, we simplify our definition of the mode to the real value $\vartheta := \arg \max_{x\in D} f(x)$; this implies the necessity to modify the previous definition of the kernel estimator of the mode by $ \vartheta_n := \arg \max_{x\in D}\hat{f}_n(x)$. 
\subsection{Some assumptions}
Now, some assumptions needed to study the asymptotic properties of the estimator $\hat{f}_n(x)$, are introduced and gathered below for easy reference.
\begin{itemize}
 \item[{A1.}] $\{X_n,n\geqslant1\}$ is a sequence of stationary WOD random variables with dominant coefficients $g_X(n)$.  
 \item[{A2.}] $\{Y_n,n\geqslant1\}$ is a sequence of stationary WOD random variables with dominant coefficients $g_Y(n)$ which   are assumed  to be independent from the random variables of interest $\{X_n,n\geqslant1\}$. 
 \item [A3.] $\int_{a_F}^{\infty}\frac{dF(z)}{G^2(z)}<+\infty. $
 \item [A4.] $G(\cdot)$ is a Lipschitz function.
 \item [A5.]  $K$ is a Lipschitz continuous probability density function satisfying:\\ $\int x K(x)dx=0$, $\int x^2K(x)dx<\infty$, $\int |x|K^2(x)dx<\infty$  and $\sup(K)<\infty$. 
\item [A6.] The sequence $h_n$ is positive and satisfies $h_n\rightarrow 0$ and $(nh_n^4)^{-1/2}\log^{1/2}(ng(n)) \rightarrow 0$ as $n\rightarrow \infty$.
\item [A7.] $f$ is twice continuously differentiable on $D$ with second derivative $f^{(2)}(\vartheta) \ne 0$ and  $\sup_{x\in D}\vert f^{(k)}(x) \vert < \infty$ for $k=1,2$.
\item [A8.] The unique mode $\vartheta$ satisfies, for any $\epsilon > 0$ and $x$, there exists a $\varrho > 0$ such that $\vert \vartheta - x\vert \geq \epsilon$ implies that $\vert f(\vartheta) - f(x)\vert \geq \varrho$.
\end{itemize} 
\begin{remark}(Discussion on the assumptions)\\
Assumptions (A1)-(A3) imply those of El Alem et al. (\cite{MZT}, 2024) and are necessary to use their results. Assumption (A4) is primarily technical and is involved in computing the fluctuation term; it has been utilized by Gheliem and Guessoum (\cite{Gheliem}, 2022) for left-truncated and associated data.    Additionally, the conditions in (A5) are fundamental requirements for the kernel function, which are satisfied by the majority of kernels, such as the Epanechnikov and Gaussian kernels. (A6) is a standard condition in nonparametric estimation of the bandwidth. Moreover, (A7) is a technical regularity condition for the density function $f(\cdot)$. Finally, assumption (A8) stipulates the uniform uniqueness of the mode point.
\end{remark}
\subsection{Strong Uniform Consistency}
Now, our first result is the strong uniform consistency  with a rate of the kernel density estimator $\hat{f}_n(x)$. To study the asymptotic behaviour of $\hat{f}_n(\cdot)$ we first study the fluctuation term $f_n(\cdot)-\mathbf{E}f_n(\cdot)$.
\begin{proposition}\label{prop1}
Under Assumptions  (A1)-(A7), we have
$$
\sup_{x\in D}\left\vert f_n(x)-\mathbf{E}\left(f_n(x)\right) \right\vert = O\left(  (nh_n^4)^{-1/2}\log^{1/2}(ng(n)) \right)\; \mathrm{a.s.}, \; \mathrm{as} \; n\rightarrow \infty,
$$
where $g(n)=max(g_X(n),g_Y(n))$.
\end{proposition}
\begin{theorem}\label{th1}
Under assumptions (A1)-(A7), we have
$$
\sup_{x\in D}\left\vert \hat{f}_n(x)-f(x) \right\vert = O\left(  (nh_n^4)^{-1/2}\log^{1/2}(ng(n)) + h_n^2  \right)\; \mathrm{a.s.}, \; \mathrm{as} \; n\rightarrow \infty.
$$
\end{theorem}
\begin{remark}
The convergence rate in this case is not as good as that reported by Wang et al. (\cite{WANG}, 2022) in the complete and WOD case, primarily due to the truncation effect. Notably, NA, NOD, NSD, and END sequences imply the WOD sequence, but the reverse does not hold true. These dependency modes represent specific cases within the scope of those studied in this paper. In instances where $g(n) = O\left(n^{\kappa}\right)$ for any $\kappa \geq 0$, the convergence rate is approximately $O\left((nh_n^4)^{-1/2}\log^{1/2}(n) + h_n^2\right)$.
\end{remark}
\subsection{Application to the mode estimate}
As an application of Theorem \ref{th1} we obtain the almost sure convergence rate of $\vartheta_n$.
\begin{theorem}\label{th2}
Under assumptions (A1)-(A9), we have
$$
\vert \vartheta_n - \vartheta \vert = O\left( (nh_n^4)^{-1/4}\log^{1/4}(ng(n)) + h_n  \right)\; \mathrm{a.s.}, \; \mathrm{as} \; n\rightarrow \infty.
$$
\end{theorem}
\begin{remark}
The estimate $\vartheta_n$ is not necessarily unique; therefore, all the results in this paper will pertain to any sequence of random variables $\vartheta_n$ satisfying $f(\vartheta_n) = \sup_{x\in \mathbb{R}} \hat{f}_n(x)$. We note that we can specify our choice by taking $\vartheta_n = \inf \left\{z\in \mathbb{R} : \hat{f}_n(z)=\sup_{x\in \mathbb{R}} \hat{f}_n(x)  \right\}$. For further details, refer to the works of Ould Sa\"{\i}d and Tatachak in their articles (\cite{OT1} and \cite{OT}, 2009).
\end{remark}
\section{Numerical illustrations}
\subsection{Truncated WOD Sample Construction}
To simulate a WOD sequence, let us consider the following real case. The first-order moving average (MA) process $X_t = Z_t -\nu Z_{t-1}$  is the result of a comprehensive study on annual temperature values measured in Basel from 1755 to 1957 ($Z_t$ is a Gaussian white noise process with mean $\mu_Z$ and standard deviation $\sigma_Z=0.7$). This time series has been examined by \cite{KAJB} and more recently by \cite{MZT}. Their study revealed that the optimal MA(1) model is the one with $\nu =0.9$, $\mu_Z=0$, and $\sigma_Z=0.7$ based on the Akaike Information Criterion (AIC). Therefore, as a linear combination of multivariate normal variables remains multivariate normal, the actual data $X_t$ possesses a multivariate normal distribution with a zero mean vector and covariance matrix $\Sigma_\nu$.
$$
\Sigma_\nu=
\begin{pmatrix}
(1+\nu^2)\sigma_Z^2 & -\nu\sigma_Z^2  & 0 & \cdots & 0 & 0 & 0 \\
-\nu\sigma_Z^2 & (1+\nu^2)\sigma_Z^2 & -\nu\sigma_Z^2 &  \cdots & 0 & 0 & 0  \\
0 & -\nu\sigma_Z^2 & (1+\nu^2)\sigma_Z^2 & \cdots & 0 &  0 & 0 \\
\vdots & \vdots & \vdots & \ddots  & \vdots & \vdots & \vdots\\
0 & 0 & 0 &  \cdots & (1+\nu^2)\sigma_Z^2 & -\nu\sigma_Z^2 &  0\\
0 & 0 & 0 &  \cdots & -\nu\sigma_Z^2 & (1+\nu^2)\sigma_Z^2 &  -\nu\sigma_Z^2 \\
0 & 0 & 0 &  \cdots & 0 & -\nu\sigma_Z^2  & (1+\nu^2)\sigma_Z^2 
\end{pmatrix}
$$
By  Joag-Dev and Proschan (\cite{JP},1983) the multivariate normal distribution is NA if the off-diagonal elements of its covariance matrix are non-positive. Then for $\nu =0.9$, we have an NA sample, which is a special case of WOD sequence. In order to get a random left truncated WOD sequence, we generate the data as follows. 
\begin{itemize}
\item[$\bullet$] \textbf{Step 1.}  The sequence $\{(X_i,Y_i), i =1,\cdots,n \}$ is generated as follows:
\begin{itemize}
\item \textbf{Variable of interest $X$}: We first generate $(N+1)$  independent and identically distributed (iid) random variables (rv's) $Z_i$ drawn from the normal distribution $\mathcal{N}(0, 0.46)$. Then generate the WOD sequence $\{X_i, i=1,\ldots, N\}$ by $X_i=Z_i-0.9 Z_{i-1}$.
\item \textbf{The truncated rv $Y$}: In the same way we compute the truncated rv's $\{Y_i, i=1,\ldots, N\}$, we first generate $(N+1)$ iid rv's $\tilde{Z}_t$ from $\mathcal{N}(0,\sigma_{\tilde{Z}}^2)$. Subsequently, we generate $Y_i=\tilde{Z}_i-\nu \tilde{Z}_{i-1}$  for $i=1,\ldots, N$, where  $0<\nu<1$. The two parameters $\sigma_{\tilde{Z}}^2$ and $\nu$ are adapted in order to control the rate of truncation $(1-\alpha)$.
\item \textbf{Observed data}:  Finally we keep the $n$ observations $\{(X_i,Y_i), i =1,\cdots,n \}$ of the couple of rv's $(X_i,Y_i)$ satisfying the condition $X_i \geq Y_i$.
\end{itemize}
\item[$\bullet$] \textbf{Step 2.} Using the simulated observed data $\{(X_i,Y_i), i =1,\cdots,n \}$, compute the kernel density estimator $\hat{f}_n(x)$.
\item[$\bullet$] \textbf{Step 3.} We repeat $M$ simulation runs as described in Steps 1 and 2 for every fixed combination of size $n$ and truncating rate $(1-\alpha)$.
\item[$\bullet$] \textbf{Step 4.} We compute  the global mean squared error (GMSE) across $M$ Monte Carlo trials, defined as 
$$
GMSE=\dfrac{1}{M H} \sum\limits_{k=1}^{M} \sum\limits_{j=1}^{H} \left( \hat{f}_{n,k} (x_{j})-f(x_{j})\right) ^{2},
$$ 
where $H$ is the number of equidistant points $ x_{j} $ belonging to the range $ [-3,3] $ and $\hat{f}_{n,k} (x_{j})$ is the value of $ \hat{f}_{n} (x_{j}) $ computed at iteration $k$.
\end{itemize}
For $n\in \{50; 100; 500\}$, the GMSE's values and the curves for the Kernel density estimator with $(1 - \alpha) \in \{10\%; 30\%; 50\%\}$ are shown in the following subsection.
\subsection{Simulation results}
\begin{figure}[H]
\begin{center}
\includegraphics[width=4.78cm]{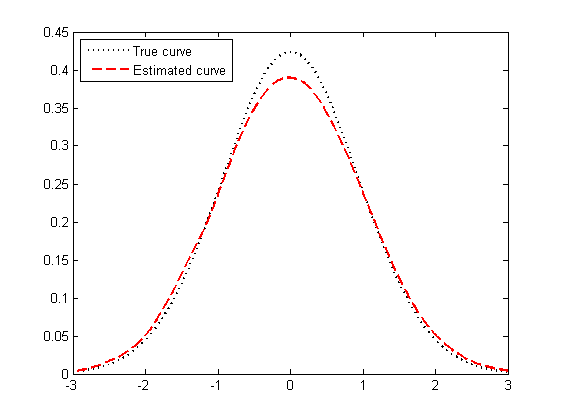}
\hspace {-.65cm}
\includegraphics[width=4.78cm]{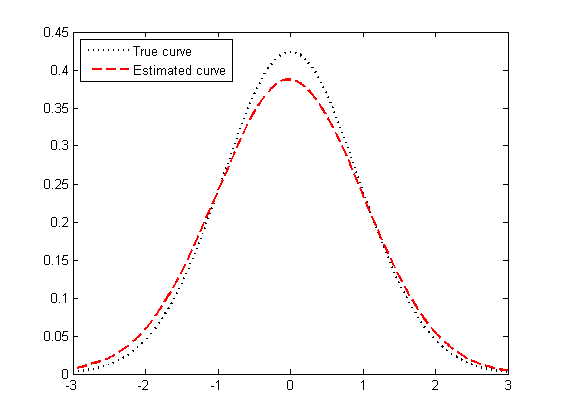}
\hspace {-.65cm}
\includegraphics[width=4.78cm]{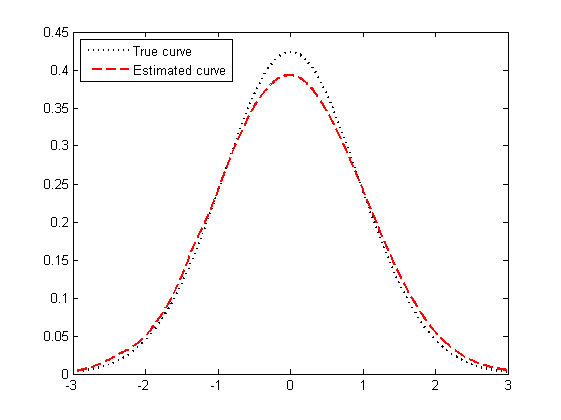}
 \caption{Kernel density estimator : n=50 and TR $\approx$ $10 \%$, $30 \%$, $50 \%$ }
 \label{fig1}
 \end{center}
\end{figure}

\begin{figure}[H]
\begin{center}
\includegraphics[width=4.78cm]{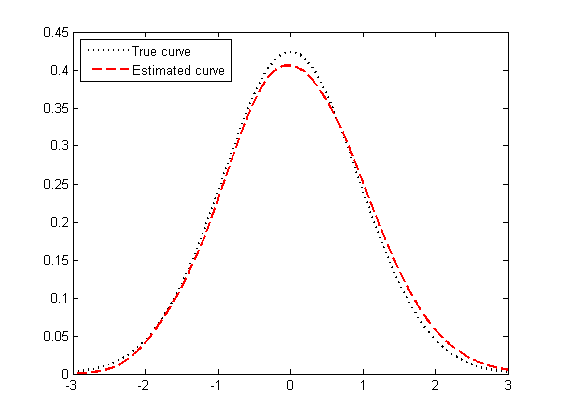}
\hspace {-.65cm}
 \includegraphics[width=4.78cm]{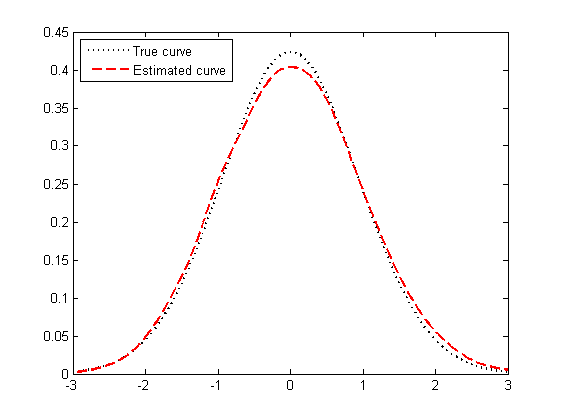}
\hspace {-.65cm}
 \includegraphics[width=4.78cm]{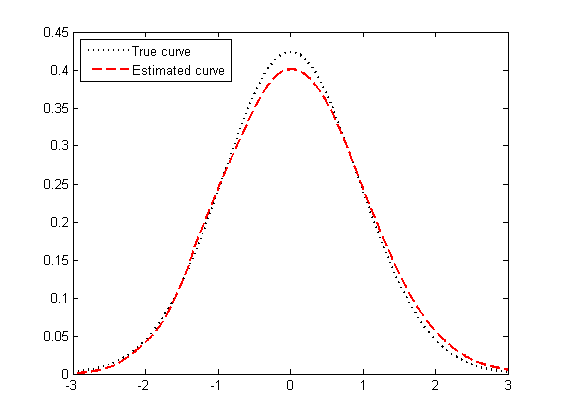}
 \caption{Kernel density estimator : n=100 and TR $\approx$ $10 \%$, $30 \%$, $50 \%$ }
 \label{fig2}
 \end{center}
\end{figure}

\begin{figure}[H]
\begin{center}
 \captionsetup{justification=centering}
\includegraphics[width=4.78cm]{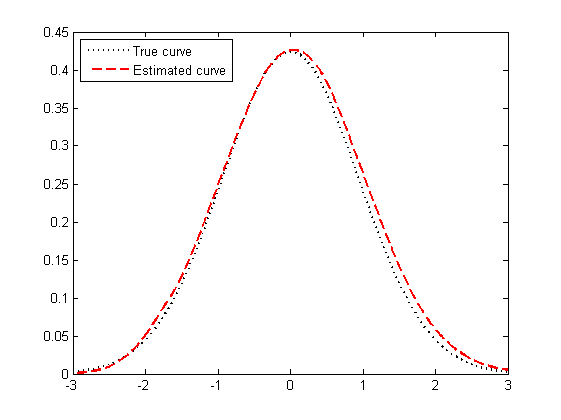}
\hspace {-.65cm}
 \includegraphics[width=4.78cm]{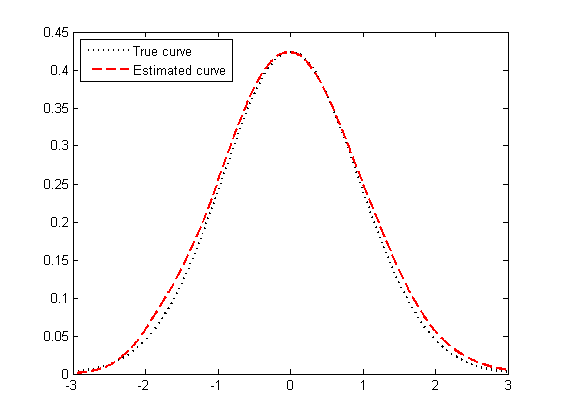}
\hspace {-.65cm}
 \includegraphics[width=4.78cm]{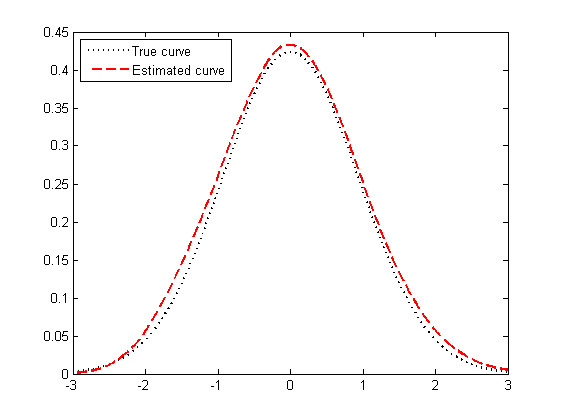}
 \caption{Kernel density estimator : n=500 and TR $\approx$ $10 \%$, $30 \%$, $50 \%$ }
 \label{fig3}
 \end{center}
\end{figure}

In this part, we examine with simulated data  the finite sample performance of our estimators $\hat{f}_n(x)$ and $\vartheta_n$ by considering some fixed-size particular cases and varying the rate of truncation.  We calculate our estimator based on the observed data $\{(X_i, Y_i), i = 1,\cdots, n\}$, by choosing a Gaussian kernel $K(x)= (1/\sqrt{2\pi}) e^{-(1/2)x}$. In all cases, following Wang et al. (\cite{WANG}, 2022), we took $h_n=O(n^{-1/5})$ for a suitable positive constant $C$, which is the optimal choice in the case of complete and WOD data, and it satisfies our standard conditions.  The GMSE is calculated for $\hat{f}_n(x)$ and the corresponding values are shown in Table \ref{table1}. In addition, to visualize how the estimator $\hat{f}_n(x)$ fits, we plot the true and the estimated curves in Figures \ref{fig1}, \ref{fig2} and \ref{fig3}.
\begin{table}[h]
\caption{GMSE values for Kernel density estimator }\label{table1}
\begin{tabular*}{\textwidth}{@{\extracolsep\fill}lccc}
\toprule%
& \multicolumn{3}{@{}c@{}}{GMSE}  \\\cmidrule{2-4}
$1-\alpha$ & $n=50$ & $n=100$ & $n=500$   \\
\midrule
10\%  & 0.0061 & 0.0031 & 0.0011  \\
30\% & 0.0073 & 0.0034 & 0.0013  \\
50\% & 0.0084 & 0.0047 & 0.0016  \\
\botrule
\end{tabular*}
\end{table}

\begin{figure}[H]
\begin{center}
\includegraphics[width=4.78cm]{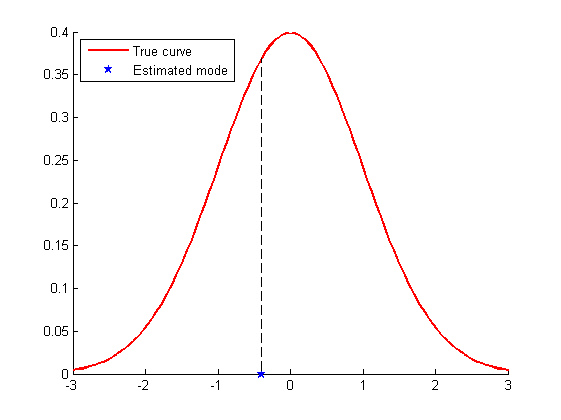}
\hspace {-.65cm}
\includegraphics[width=4.78cm]{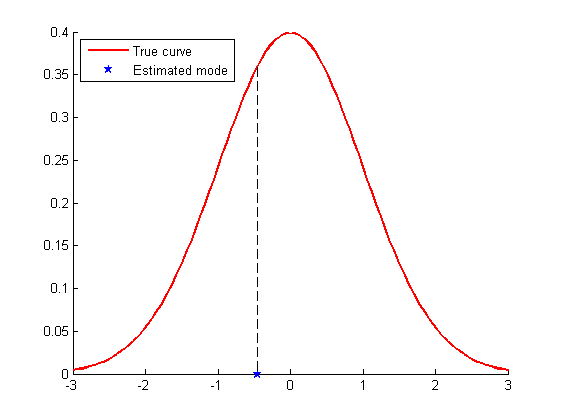}
\hspace {-.65cm}
\includegraphics[width=4.78cm]{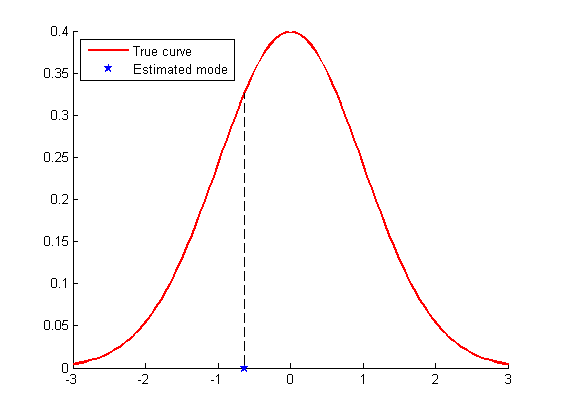}
 \caption{Kernel mode estimator : n=50 and TR $\approx$ $10 \%$, $30 \%$, $50 \%$ }
 \label{fig11}
 \end{center}
\end{figure}

\begin{figure}[H]
\begin{center}
\includegraphics[width=4.78cm]{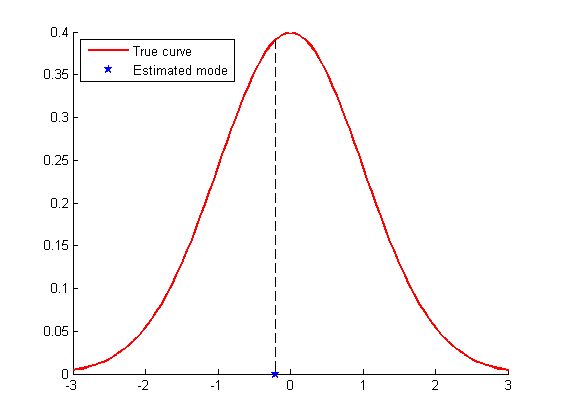}
\hspace {-.65cm}
 \includegraphics[width=4.78cm]{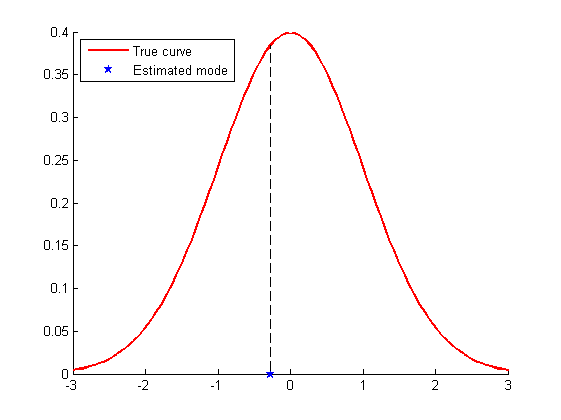}
\hspace {-.65cm}
 \includegraphics[width=4.78cm]{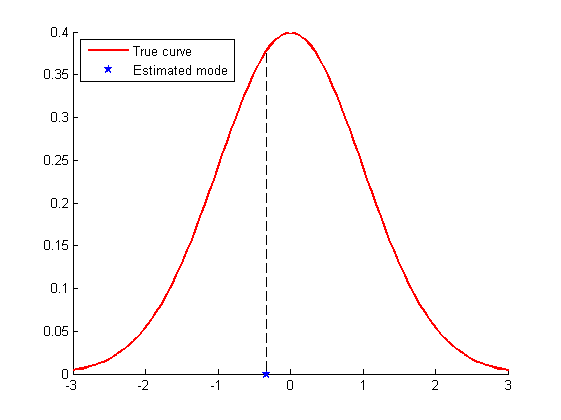}
 \caption{Kernel mode estimator : n=100 and TR $\approx$ $10 \%$, $30 \%$, $50 \%$ }
 \label{fig12}
 \end{center}
\end{figure}

\begin{figure}[H]
\begin{center}
 \captionsetup{justification=centering}
\includegraphics[width=4.78cm]{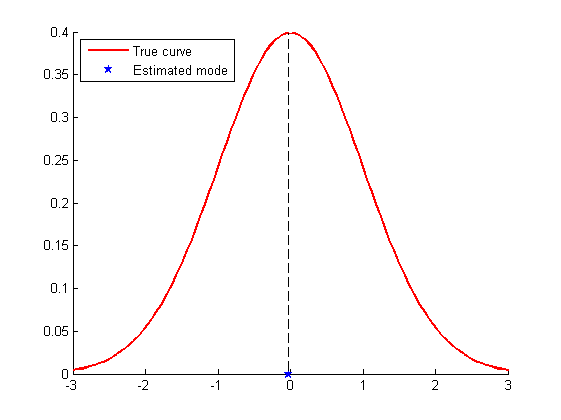}
\hspace {-.65cm}
 \includegraphics[width=4.78cm]{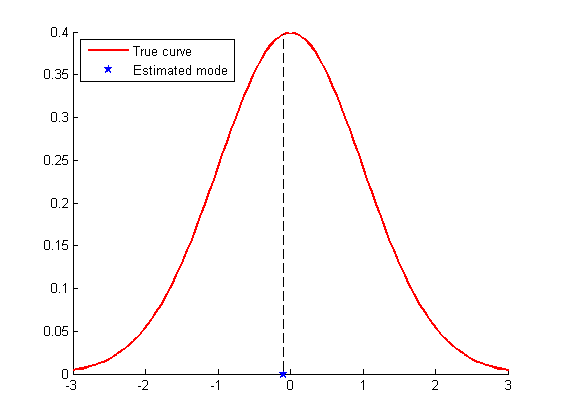}
\hspace {-.65cm}
 \includegraphics[width=4.78cm]{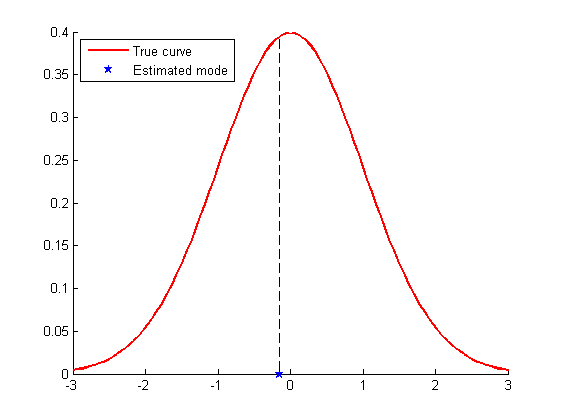}
 \caption{Kernel mode estimator : n=500 and TR $\approx$ $10 \%$, $30 \%$, $50 \%$ }
 \label{fig13}
 \end{center}
\end{figure}

To underscore the theoretical result presented in Theorem \ref{th2}, we perform a simulation study to examine the behavior of the kernel estimator for the mode. The variables \( X \) and \( Y \) are generated using the same procedure as described previously. This simulation is repeated \( M = 300 \) times, with different sample sizes \( n \) and varying truncation rates \( 1 - \alpha \), as detailed in Table \ref{table2} and depicted in Figures \ref{fig11}, \ref{fig12}, and \ref{fig13}. The table reports the mean squared error (MSE) of the mode estimator. The findings reveal that as the truncation rate increases, the MSE worsens, whereas the MSE improves with larger sample sizes.
\begin{table}[h]
\caption{MSE values for Kernel mode estimator }\label{table2}
\begin{tabular*}{\textwidth}{@{\extracolsep\fill}lccc}
\toprule%
& \multicolumn{3}{@{}c@{}}{MSE}  \\\cmidrule{2-4}
$1-\alpha$ & $n=50$ & $n=100$ & $n=500$   \\
\midrule
10\%  & $4.51\times 10^{-5}$ & $2.45\times 10^{-5}$ & $1.23\times 10^{-5}$  \\
30\% & $6.27\times 10^{-5}$ & $2.78\times 10^{-5}$ & $1.47\times 10^{-5}$  \\
50\% & $8.33\times 10^{-5}$ & $3.92\times 10^{-5}$ & $1.89\times 10^{-5}$  \\
\botrule
\end{tabular*}
\end{table}

\subsection{Discussion}
\begin{itemize}
\item We can see through the simulation results that the quality of fit to the theoretical curve deteriorates slightly when the truncation rate $1-\alpha$ increases.
\item We can thus observe that the quality of fit improves when $n$ increases, ie when $n$ tends towards infinity, a better fit is obtained.
\item The simulations revealed that the estimator is less affected by the truncation rate $1-\alpha$ than by the small sample size.
\end{itemize}
\section{Analysis of real data}
To highlight the pertinence of this study through a real-life example, we apply our proposed estimator to the car brake pedal lifetime data, which can be found in Lawless (\cite{Lawless}, 2002). This data set  include the brake pads lifetimes of 98 individual cars. It can be deduced that the pads which possess longer lifetimes own greater chance of being observed in the selected sample. This bias happens owing to non-random sampling of components and therefore the presence of left truncation variable. Thus, the sample inclusion criterion is $X\geq Y$, where $Y$ is the number of kilometers driven for brake pads at the time of sampling, and $X$ is the number of kilometers driven for brake pads at the time of failure. The aforementioned brake pads data are written as $\{(X_i, Y_i)\;, X_i \geq Y_i, i = 1, 2,\cdots , n\}$ with $n = 98$. By utilizing MATLAB software, we have determined that the given dataset of $X_i$ conforms to the gamma distribution through the Kolmogorov-Smirnov (K-S) test. In this context, the null hypothesis posits that the dataset is derived from a gamma distribution, while the alternative hypothesis suggests otherwise. Upon analyzing the data regarding the lifespan of car brake pads, the K-S test statistic is $0.038141$, with a corresponding p-value of $0.99387$. In this instance, the K-S critical value of $0.13738$ for a $5\%$ level of significance surpasses the calculated statistic, and the level of significance is smaller compared to the p-value. Consequently, the affirmation of the null hypothesis implies that the gamma distribution constitutes a plausible fit for this dataset. To ascertain the parameters of the gamma distribution based on the provided data, we employ maximum likelihood estimation (MLE), resulting in $({\delta}_1=6.5768,\; {\beta}_1= 10.2982)$. Similarly, we prove that the data set $Y_i$ conforms to the gamma distribution, with parameters $({\delta}_2=6,9,\; {\beta}_2= 3.3308)$, using the Kolmogorov-Smirnov (K-S) test. The K-S test statistic is $0.079748$, with a corresponding p-value of $0.53493$. In this case, the critical K-S value of $0.13738$ for a significance level of $5\%$ exceeds the calculated statistic, and the significance level is smaller than the p-value.
The results of the chi-squared test indicate a p-value of $0.1834$. Since this p-value is above the $5\%$ significance level, it suggests that there is insufficient evidence to reject the null hypothesis of independence between the $X_i$ and $Y_i$ variables.
\begin{figure}[H]
\begin{center}
 \captionsetup{justification=centering}
\includegraphics[width=5.78cm]{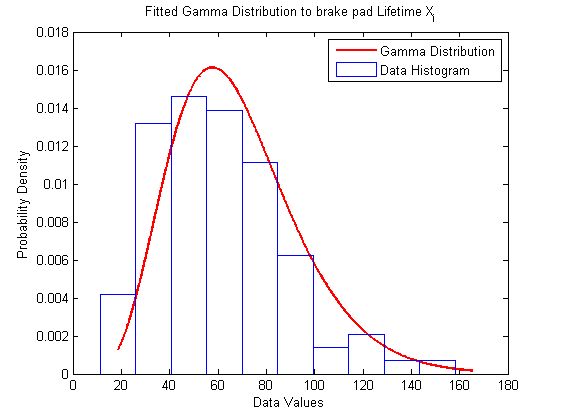}
\includegraphics[width=5.78cm]{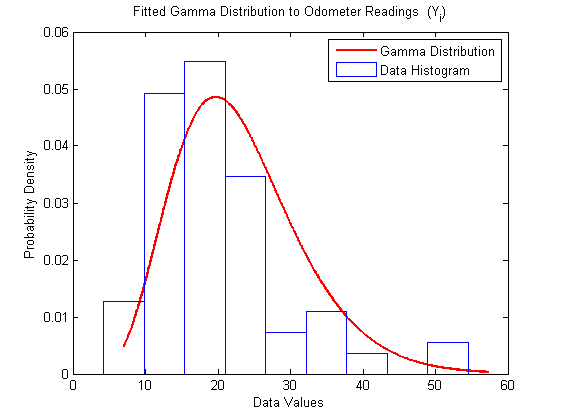}
 \caption{Fitted Gamma Distribution to Brake Pad Lifetime $X_i$ and Odometer Readings of 98 automobiles.}\label{figG}
 \end{center}
\end{figure}
\noindent
The proposed estimator $\hat{f}_n$ can be applied to this dataset with $n=98$ and truncation rate $1-\alpha \approx 10 \%, 30 \%, 50 \%$ by generating, as discussed before, the $X_i$ according to the gamma$({\alpha}_1,\; {\beta}_1)$ distribution and the $Y_i$ according also to a gamma distribution where its parameters vary in order to control the truncation rate. Subsequently, we retain the data $(X_i, Y_i)$ such that $X_i\geq Y_i$. It is important to note that these samples satisfy the conditions of our model since the WOD variables contain the independent variables as a special case and $X_i$ is independent of $Y_i$, for example, we are referring to (\cite{KM}, 2021).  We take the bandwidth parameter $h_n = O\left(n^{-1/5} \right)$, which satisfies condition (H), and the Gaussian kernel function is used.
\begin{figure}[H]
\begin{center}
 \captionsetup{justification=centering}
\includegraphics[width=4.78cm]{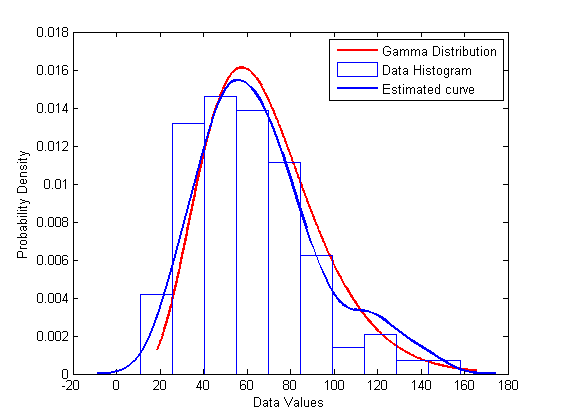}
\includegraphics[width=4.78cm]{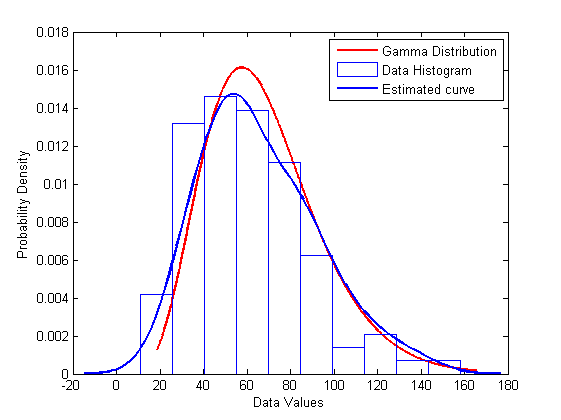}
\includegraphics[width=4.78cm]{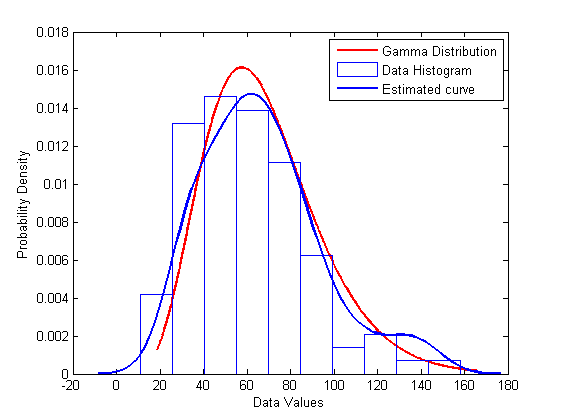}
 \caption{Kernel density estimator for Lifetime of automobile brake pads :  $1-\alpha \approx 10 \%, 30 \%, 50 \%$, respectively.}
 \label{figGM}
 \end{center}
\end{figure}
\noindent
Figure \ref{figGM} clearly shows the good fit of our estimator to the fitted parametric distribution using the maximum likelihood method. On the other hand, we can also conclude that our nonparametric estimator fits the histogram of the actual data better than the parametric estimator, and this for the different truncation rates.
\section{Proofs and lemmas}
\textbf{Proof of Proposition~\ref{prop1}.}
The idea consists in using an exponential inequality taking into account the WOD structure. The compact set $D$  can be covered by a finite number $q_n$ of intervals of length $w_n={(n^{-1}h_n^3)}^{1/2}$. Let $I_k=I(x_k,w_n)$; $k=1,\cdots , q_n$, denote each interval centered at some point $x_k$. Since $D$ is bounded, there exists a constant $M$ such that $|x-x_k|\leq w_n$. We start by writing 
\begin{equation*}\label{inegal}
\Delta_i(x)=\frac{\alpha h_n}{G(X_i)}K\left( \frac{x-X_i}{h_n} \right)-\mathbf{E}\left( \frac{\alpha h_n}{G(X_i)}K\left( \frac{x-X_i}{h_n} \right) \right)
\end{equation*}
and consider the following decomposition 
\begin{eqnarray}
\sup_{x\in D}\left\vert f_n(x)-\mathbf{E}\left(f_n(x)\right) \right\vert &= & \sup_{x\in D}\left\vert  \frac{1}{nh_n^2}\sum_{i=1}^n\left( \Delta_i(x) - \Delta_i(x_k)\right) + \frac{1}{nh_n^2}\sum_{i=1}^n\Delta_i(x_k)  \right\vert \notag  \\
            &\leq & \max_{1\leq k\leq q_n}\sup_{x\in D}\frac{1}{nh_n^2}\sum_{i=1}^n\left\vert \Delta_i(x) - \Delta_i(x_k) \right\vert + \max_{1\leq k\leq q_n} \frac{1}{nh_n^2}\left\vert \sum_{i=1}^n\Delta_i(x_k)  \right\vert  \notag \\
            &= & J_{1n}+J_{2n}. \notag
\end{eqnarray}
To treat $J_{1n}$, according to conditions A5 and A6, we write
\begin{eqnarray}
\frac{1}{nh_n^2}\sum_{i=1}^n\left\vert \Delta_i(x) - \Delta_i(x_k) \right\vert &\leq & \frac{1}{nh_n} \sum_{i=1}^n \frac{\alpha}{G(X_i)}\left\vert K\left( \frac{x-X_i}{h_n} \right)-K\left( \frac{x_k-X_i}{h_n} \right)\right\vert\notag  \\
            &\qquad +& \frac{\alpha}{h_n}\mathbf{E}\left( \frac{1}{G(X_i)}\left\vert K\left( \frac{x-X_i}{h_n} \right)-K\left( \frac{x_k-X_i}{h_n} \right)\right\vert \right)   \notag \\
            &\leq & \frac{C\alpha \vert x-x_k \vert}{h_n^2a_F} \notag \\
            &\leq & \frac{C\alpha w_n}{h_n^{2}}  = O\left(\frac{1}{\sqrt{nh_n}}\right), \notag
\end{eqnarray}
so, $J_{1n}= O\left( \frac{1}{\sqrt{nh_n}} \right)$.\\
Now to deal with $J_{2n}$, set $L_x(y) := \frac{\alpha h_nK\left((x-y)/h_n\right)}{G(y)}$, 
and note that as $K$ and $G$ are $M_1$-lipschitzian and $M_2$-lipschitzian respectively, with $K$ upper bounded (from Assumption A5) and $G(y) \geq a_F >0 $, then for any $z, y \in D$, and for a given $h_n$, we have
\begin{eqnarray*}
\left|L_x(z)-L_x(y)\right|&=&\alpha h_n\left|\frac{K((x-z)/h_n)}{G(z)}-\frac{K((x-y)/h_n)}{G(y)}\right| \notag \\
&\leq& \left(\frac{\alpha M_1}{a_F} + \alpha h_nM_2\frac{\left\|K\right\|_{\infty}}{a_F^2}\right)\left|z-y\right|,\quad \text{where}\; \left\|K\right\|_{\infty}=\sup(K).
\end{eqnarray*}
As $h_n\rightarrow 0$, we can find a positive constant $M^{'}$ such that $\left(\frac{M_1}{a_F} + h_nM_2\frac{\left\|K\right\|_{\infty}}{a_F^2}\right)\leq M^{'}$. This inequality implies the Lipschitz continuity of $L_x(\cdot)$ on $D$, and consequently, it is of bounded variation. Thus, there exist monotonic functions $L_{1}(\cdot)$ and $L_{2}(\cdot)$ such that 
\begin{equation}\label{eqL}
L_{x}(y)=L_{1}(y)-L_{2}(y).
\end{equation}
We now make use of the following exponential-type inequality.
\begin{lemma}\label{lem5} (Theorem 2.1 in Liu et al. (\cite{Liu}, 2016)). 
Let $\{X_n,\; n\geqslant 1\}$ be a sequence of identically distributed WOD rv's with $\mathbb{E}(X_n)= 0$ and  $\vert X_n \vert \leq c$ for each $n\geq 1$, where $c$ is a positive constant.
Then for any $ 0 < B < 1 $  and  $ 0 < \epsilon < cB/(1-B)$, we have 
$$
\mathbb{P}\left(\frac{1}{n}\left\vert \sum_{i=1}^n  X_i \right\vert >\epsilon \right)\leq 2g(n)\exp\left\lbrace -\frac{n\epsilon^2}{4C} \right\rbrace,
$$
where $C=\frac{c^2}{2(1-B)}$.
\end{lemma}
\begin{remark}
The use of this Lemma is motivated by its simple computational aspect compared to alternatives such those of Xia et al. (\cite{Xia}, 2018) and Shen (\cite{Shen}, 2013).
\end{remark}
According to \ref{eqL},  we require the following notations.\\
\begin{equation*}
L(X_i, x_k) :=L_{x_k}(X_i) \,\text{ and } f_{p,n}(x_k) := \frac{1}{nh_n^2}\sum_{i=1}^nL_p\left(X_i, x_k \right);\; p=1,2.
\end{equation*}
So, it is easily seen that $f_{n}(x_k)=f_{1,n}(x_k)-f_{2,n}(x_k)$.\\
Then for each $n \geq 1$; $i = 1,\ldots ,n$ and  $p = 1,2$, let
$$
\Delta_{p,i}(x_k)=L_p\left( X_i, x_k \right)-\mathbf{E}\left( L_p\left( X_i, x_k \right) \right).
$$
It follows that
$$
f_{p,n}(x_k)-\mathbf{E}\left(f_{p,n}(x_k)\right)=\frac{1}{nh_n^2}\sum_{i=1}^n \Delta_{p,i}(x_k)$$
and
$$
\frac{1}{nh_n^2}\sum_{i=1}^n \Delta_{i}(x_k)= \frac{1}{nh_n^2}\sum_{i=1}^n \Delta_{1,i}(x_k)- \frac{1}{nh_n^2}\sum_{i=1}^n \Delta_{2,i}(x_k).
$$
Remark that
\begin{eqnarray}
\frac{1}{nh_n^2}\left|\sum_{i=1}^n \Delta_{i}(x_k)\right| &\leq & \frac{1}{nh_n^2}\left|\sum_{i=1}^n \Delta_{1,i}(x_k)\right| + \frac{1}{nh_n^2}\left|\sum_{i=1}^n \Delta_{2,i}(x_k)\right|.
\end{eqnarray}
Recall that if $\{X_1, \ldots, X_n\}$ is WOD, then $\{{\Delta}_{p,1}(x_k), \ldots, {\Delta}_{p,n}(x_k)\}$  is also WOD. Notably, $\mathbf{E}(\Delta_{p,i}(x_k))=0$.  According to Kolomogorov and Formin (\cite{kol}, 1975), we can adopt $L_{1}(y)=V_a^y(L_x)$ where $V_a^y(L_x)$ represents the bounded variation of $L_x$ on the compact $[a, y]$ (with $a\leq y\leq b$, where $a$ and $b$ are the end-points of the process $\{X_i\}$, denoted as $a_F$ and $b_F$).\\
According to the Corollary on p330 in \cite{kol}, the function $v(y)=V_a^y$ is increasing. Hence, $V_a^y(L_x) \leq V_a^b(L_x)\leq c_1$, which implies $L_{1}(y)\leq c_1$ based on their definitions 1 and 2 on p328. Assuming that $L_{x}(y)$ is bounded (i.e., there exists $c_2$ such that $\left\vert L_{1}(y)\right\vert \leqslant c_2$), we can establish $\left\vert L_{2}(y)\right\vert \leq \left\vert L_{1}(y) - L_x(y)\right\vert \leq \left\vert L_{1}(y)\right\vert + \left\vert L_{x}(y)\right\vert \leq c_1+c_2$. Consequently, $\vert {\Delta}_{p,i}(x_k) \vert \leq c = 2(c_1 + c_2)$. 
This implies that $\left\lbrace {\Delta}_{p,i}(x_k), i=1,\ldots, n\right\rbrace$ satisfies the conditions of Lemma \ref{lem5} and then we have
\begin{eqnarray}\label{ineg}
\mathbf{P}\left(\frac{1}{n}\left\vert \sum_{i=1}^n \Delta_i(x_k) \right\vert >\epsilon h_n^2 \right) &\leq & \mathbf{P}\left(\frac{1}{n} \left\vert \sum_{i=1}^n {\Delta}_{1,i}(x_k) \right\vert >\frac{\epsilon h_n^2}{2} \right)+\mathbf{P}\left(\frac{1}{n} \left\vert \sum_{i=1}^n  {\Delta}_{2,i}(x_k) \right\vert >\frac{\epsilon h_n^2}{2} \right) \notag  \\
            &\leq & 4g(n)\exp\left\lbrace -\frac{n\epsilon^2 h_n^4}{16C} \right\rbrace. 
\end{eqnarray}
Next, if we choose $\epsilon=\epsilon_o\sqrt{\frac{\log(ng(n))}{nh_n^4}}$ for all $\epsilon_o>0$ and applying (\ref{ineg}), we write 
\begin{eqnarray}\label{borel}
\mathbf{P}\left(\max_{1\leq k\leq q_n}\frac{1}{n}\left\vert \sum_{i=1}^n \Delta_i(x_k)\right\vert >h_n^2\epsilon_o\sqrt{\frac{\log(ng(n))}{nh_n^4}}\right) &\leq & \sum_{i=1}^{q_n} \mathbf{P}\left(\left\vert \sum_{i=1}^n \Delta_i(x_k)\right\vert  >nh_n^2\epsilon_o\sqrt{\frac{\log(ng(n))}{nh_n^4}}\right) \notag  \\
            &\leq & \sum_{i=1}^{q_n}4g(n)\exp\left\lbrace -\frac{nh_n^4{\epsilon_o}^2\dfrac{\log(ng(n))}{nh_n^4}}{16C} \right\rbrace \notag \\
            &\leq & 4q_ng(n)\exp\left\lbrace -\frac{{\epsilon_o}^2\log(ng(n))}{16C} \right\rbrace \notag \\
            &\leq & 4M{(w_n)}^{-1}g(n){(ng(n))}^{-c_o{\epsilon_0}^2}
 \notag \\
            &\leq & \frac{4M}{\sqrt{{(nh_n)}^3}}n^{-c_o{\epsilon_o}^2 +2}.
\end{eqnarray}
By Assumption (A2) and for a suitable choice of  $\epsilon_o$ (i.e ${\epsilon_o}^2>\frac{3}{c_o}$), the right hand side term in (\ref{borel}) is the general term of a convergent series. Then by the Borel-Cantelli's lemma we get
$$
J_{2n}=O\left(\sqrt{\frac{\log(ng(n))}{nh_n^4}} \right)\; \text{a.s., as} \; n\rightarrow \infty,$$
which ends the proof.  $\Box$\\
\textbf{Proof of Theorem~\ref{th1}.}
The proof is based on triangular inequality hereafter, and is broken into proofs of the following lemmas.\\
First observe that 
\begin{eqnarray}
|\hat{f}_n(x) - f(x)| &\leq & |\mathbf{E}[f_n(x)] - f(x)| + |\hat{f}_n(x) - f_n(x)| + |f_n(x) - \mathbf{E}[f_n(x)]|  \notag \\                   
              &=& R_{1n} \quad + \quad R_{2n} \quad + \quad  R_{3n}. \notag 
\end{eqnarray}
\begin{lemma}\label{lem1}
Assume that hypotheses (A5) and (A7) hold, then
$$
\sup_{x\in D}|\mathbf{E}[f_n(x)] - f(x)|=O\left(h_n^2 \right)
$$
\end{lemma}
\textbf{Proof of Lemma~\ref{lem1}.}
The asymptotic behavior of $R_{1n}$ is standard, in the sense that it is not affected by the dependence structure. Indeed, using a change of variable and a Taylor expansion, we have
\begin{eqnarray}
\mathbf{E}[f_n(x)] - f(x) & = &\frac{1}{h_n}\int \frac{\alpha}{G(t)}K\left(\frac{x-t}{h_n}\right)d\mathbf{F}(t) - f(x)  \notag \\                   
              &=& \int K(z)\frac{{(zh_n)}^2}{2}f^{(2)}(x^*)dz , \notag 
\end{eqnarray}
with $x^*$ is between $x-zh_n$ and $x$. Thus
$$
|\mathbf{E}[f_n(x)] - f(x)|\leq \frac{{h_n}^2}{2}\sup_{x\in D}|f^{(2)}(x)|\int z^2K(z)dz.
$$
Under the given conditions, the result holds. $\Box$\\
For the term $R_{2n}$ and $R_{3n}$, we have the following results :
\begin{lemma}\label{lem2}
Assume that Assumptions (A1)-(A3)  hold, then
$$
\sup_{x\in D}|\hat{f}_n(x)-f_n(x)|= O\left(\sqrt{\frac{\log(ng(n))}{n}}\right)
$$
\end{lemma}
\textbf{Proof of Lemma~\ref{lem2}.}
We write
\begin{eqnarray}
|\hat{f}_n(x)-f_n(x)| & = &\left| \frac{1}{nh_n}\sum_{i=1}^n \left[ \frac{\alpha_n}{G_n(X_i)}-\frac{\alpha}{G(X_i)}\right] K\left(\frac{x-X_i}{h_n}\right)  \right|  \notag \\ 
                    & \leq &\frac{|\alpha_n - \alpha|}{nh_n}\sum_{i=1}^n\frac{1}{G_n(X_i)}K\left(\frac{x-X_i}{h_n}  \right)+\frac{\alpha}{nh_n}\sum_{i=1}^nK\left(\frac{x-X_i}{h_n} \right)\left|\frac{1}{G_n(X_i)} -\frac{1}{G(X_i)}\right|  \notag \\                   
              &\leq& \left\{\frac{|\alpha_n - \alpha|}{G_n(a_F)}+\frac{\alpha}{G_n(a_F)G(a_F)}\sup_{x\in D}|G_n(x)-G(x)|  \right\}\frac{1}{nh_n}\sum_{i=1}^nK\left(\frac{x-X_i}{h_n} \right). \notag 
\end{eqnarray}
By using Markov's inequality, a change of variables and the definition of the mode and for $\epsilon > 0$ we get
\begin{eqnarray}
\mathbf{P}\left( \frac{1}{nh_n}\sum_{i=1}^nK\left(\frac{x-X_i}{h_n}\right)\geq \epsilon \right) & \leq & \frac{\mathbf{E}\left( \frac{1}{nh_n}\sum_{i=1}^nK\left(\frac{x-X_i}{h_n}\right)\right)}{\epsilon} \notag \\ 
                    & = &\frac{1}{\epsilon h_n}\mathbf{E}\left(K\left(\frac{x-X_1}{h_n}\right)\right)  \notag \\                   
              &= & \frac{1}{\epsilon h_n}\int_R K\left(\frac{x-t}{h_n}\right)f(t)dt   \notag \\
              &= & \frac{1}{\epsilon}\int_R K\left(z\right)f(x-zh_n)dz   \notag \\
              &\leq & \frac{f(\vartheta)}{\epsilon} = O(1).  \notag 
\end{eqnarray}
Then, by following El Alem,  Guessoum and Tatachak (\cite{MZT} 2024) and using assumptions  (A1)-(A3) we get
$$
|\alpha_n - \alpha| = O\left(\sqrt{\frac{\log(ng(n))}{n}}\right) \quad \text{and} \quad \sup_{x\in D}|G_n(x)-G(x)| = O\left(\sqrt{\frac{\log(ng(n))}{n}}\right).
$$
Finally, we deduce the result. $\Box$\\
The proof of Theorem \ref{th1} is derived by combining Proposition \ref{prop1} with Lemmas \ref{lem1} and \ref{lem2}. $\Box$\\
\textbf{Proof of Theorem~\ref{th2}.}
Standard arguments give us
\begin{eqnarray}\label{eqM}
\left\vert f(\vartheta_n) - f(\vartheta)\right\vert &\leq & \left\vert f(\vartheta_n) - \hat{f}_n(\vartheta_n)\right\vert + \left\vert \hat{f}_n(\vartheta_n) - f(\vartheta)\right\vert \notag \\ 
                    & \leq & \sup_{x\in D}\left\vert \hat{f}_n(x) - f(x)\right\vert + \left\vert \hat{f}_n(\vartheta_n) - f(\vartheta)\right\vert  \notag \\                   
              &\leq & 2\sup_{x\in D}\left\vert \hat{f}_n(x) - f(x)\right\vert. 
\end{eqnarray}
Now, a Taylor expansion of $f(\cdot)$ in a neighborhood of $\vartheta$ gives
$$
f(\vartheta_n) = f(\vartheta) +  \frac{1}{2} (\vartheta_n - \vartheta)^2 f^{(2)}(\tilde{\vartheta}),
$$
where $\tilde{\vartheta}$  is between $\vartheta_n$ and $\vartheta$. Then by (\ref{eqM}), (A7) and (A8), we have
$$
\left\vert \vartheta_n - \vartheta \right\vert \leq 2\sqrt{\frac{\sup_{x\in D}\left\vert \hat{f}_n(x) - f(x)\right\vert}{f^{(2)}(\tilde{\vartheta})}}
$$
Hence, by Theorem \ref{th1}, we complete the proof of Theorem \ref{th2}. $\Box$\\
\section{Conclusion}
The motivation for this article is based on the fact that this type of incomplete data (RLT) and this notion of widely orthant dependence (WOD) are often encountered with wide application in many fields. However, our main results generalize the corresponding ones for independent samples and some negatively dependent samples such that negatively associated (NA), negatively orthant dependent (NOD), extended negatively dependent (END) and superadditive negatively dependent (SND).\\

In this paper, we establish the strong uniform consistency and the convergence rate for the kernel estimator of the probability density function in the case of left truncation and widely dependent samples. As an application, we derive the strong consistency rate for mode estimation. Additionally, a numerical illustration is conducted to evaluate the performance of the kernel estimator in a finite sample. Furthermore, a real data example is considered to support the good fit of the proposed estimator to the real density. These numerical studies show that the goodness of fit improves with an increasing sample size and deteriorates with an increasing truncation rate.\\

Another direction for future research could involve investigating the asymptotic normality of the estimator under investigation. In that respect, we could examine the behavior of other density estimators, such as the recursive estimator and the adaptive estimator, and then compare them with each other, still within the framework of this model of left truncation and widely dependence.

\section*{Declarations}
\textbf{Conflict of interest} The corresponding author, on behalf of all the authors, declares that there are no conflicts of interest.


\end{document}